\documentclass{amsart} 

\newtheorem{theorem}{Theorem}[section]
\newtheorem{proposition}[theorem]{Proposition}
\newtheorem{definition}[theorem]{Definition}
\newtheorem{lemma}[theorem]{Lemma}
\newtheorem{corollary}[theorem]{Corollary}
\newtheorem{remark}[theorem]{Remark}
\newtheorem{example}[theorem]{Example}

\newcommand{\Hn}{\hat{0}}
\newcommand{\He}{\hat{1}}
\newcommand{\Z}{\mathbb{Z}}
\newcommand{\N}{\mathbb{N}}
\newcommand{\oP}{\overline{P}}
\newcommand{\Br}{\mathrm{Br}}
\newcommand{\twist}{\mathfrak{I}}
\newcommand{\id}{\mathrm{id}}
\newcommand{\s}{\underline{s}}
\newcommand{\wh}{\widehat}

\newcommand{\SSS}{\underline{S}}
\newcommand{\wk}{\preceq}
\newcommand{\wkW}{\leq_R}
\newcommand{\Wwk}{\leq_L}
\newcommand{\Wk}{\text{Wk}}
\newcommand{\des}{\text{des}}
\begin{document}

\title[Combinatorics of twisted involutions]{The combinatorics of twisted involutions in Coxeter groups}  
\author{Axel Hultman}
\thanks{Supported by the European Commission's IHRP Programme, grant HPRN-CT-2001-00272, ``Algebraic Combinatorics in Europe''.}
\address{Department of mathematics, KTH, SE-100 44 Stockholm, Sweden}
\email{axel@math.kth.se} 
\date{April 25, 2005}

\begin{abstract}
The open intervals in the Bruhat order on twisted involutions in a Coxeter group are shown to be PL spheres. This implies results conjectured by F.\ Incitti and sharpens the known fact that these posets are Gorenstein$^*$ over $\Z_2$.

We also introduce a Boolean cell complex which is an analogue for twisted involutions of the Coxeter complex. Several classical Coxeter complex properties are shared by our complex. When the group is finite, it is a shellable sphere, shelling orders being given by the linear extensions of the weak order on twisted involutions. Furthermore, the $h$-polynomial of the complex coincides with the polynomial counting twisted involutions by descents. In particular, this gives a type independent proof that the latter is symmetric.
\end{abstract} 

\maketitle

\section{Introduction}
Let $(W,S)$ be a finite Coxeter system with an involutive automorphism $\theta$. In \cite{springer}, Springer studied the combinatorics of the twisted involutions $\twist(\theta)$. Together with Richardson, he refined his results in \cite{RS, RS2} and put them to use in the study of the subposet of the Bruhat order on $W$ induced by $\twist(\theta)$. One of their tools was another partial order on $\twist(\theta)$ which they called the {\em weak order} for reasons that will be explained later. Their motivation was an intimate connection between the Bruhat order on $\twist(\theta)$ and Bruhat decompositions of certain symmetric varieties.

The purpose of the present article is to investigate the properties of the Bruhat order and the weak order on $\twist(\theta)$ in an arbitrary Coxeter system. The Bruhat order and the two-sided weak order on a Coxeter group appear as special cases of these posets, and many properties carry over from the special cases to the general situation.

Specifically, we prove for the Bruhat order on $\twist(\theta)$ that the order complex of every open interval is a PL sphere. When $\theta = \id$, this is the principal consequence of a conjecture of Incitti \cite{incitti3} predicting that the poset is EL-shellable. The result sharpens \cite[Theorem 4.2]{hultman} which asserts that every interval in the Bruhat order on $\twist(\theta)$ is Gorenstein$^*$ over $\Z_2 = \Z/2\Z$.

Regarding the weak order on $\twist(\theta)$, we construct from it a Boolean cell complex $\Delta_\theta$ analogously to how the Coxeter complex $\Delta_W$ is constructed from the weak order on $W$. Counterparts for $\Delta_\theta$ of several important properties of $\Delta_W$ are proved. When $W$ is finite, $\Delta_\theta$ is a shellable sphere; any linear extension of the weak order on $\twist(\theta)$ is a shelling order. Moreover, the $h$-polynomial of $\Delta_\theta$ turns out to coincide with the generating function counting the elements of $\twist(\theta)$ by their number of descents. This yields a uniform proof that these polynomials are symmetric.

The remainder of the paper is organized as follows. In Section \ref{se:preliminaries}, we review necessary background material from combinatorial topology and the theory of Coxeter groups. Thereafter, in Section \ref{se:comb}, we prove some facts about the combinatorics of $\twist(\theta)$ that we need in the sequel. Most of these are extensions to arbitrary Coxeter groups of results from \cite{RS, RS2,springer}. The Bruhat order on $\twist(\theta)$ is studied in Section \ref{se:bruhat}. Finally, in Section \ref{se:weak}, we focus on the weak order on $\twist(\theta)$ and the aforementioned analogue of the Coxeter complex.

\vspace{5pt}

\noindent {\bf Acknowledgement.} The author would like to thank Mark Dukes for valuable comments on an earlier version of the paper. 

\section{Preliminaries}\label{se:preliminaries}
Here, we collect some background terminology and facts for later use.

\subsection{Posets and combinatorial topology}
A poset is {\em bounded} if it has unique maximal and minimal elements, denoted $\He$ and $\Hn$, respectively. If $P$ is bounded, then its {\em proper part} is the subposet $\oP = P\setminus\{\Hn,\He\}$.

\begin{definition}
A poset $P$ is {\em Eulerian} if it is bounded, graded and finite, and its M\"obius function $\mu$ satisfies 
\[
\mu(p,q) = (-1)^{r(q)-r(p)}
\] 
for all $p\leq q\in P$, where $r$ is the rank function of $P$.
\end{definition}

To any poset $P$, we associate the {\em order complex} $\Delta(P)$. It is the simplicial complex whose simplices are the chains in $P$. When we speak of topological properties of a poset $P$, we have the properties of $\Delta(P)$ in mind.

\begin{definition}
A poset $P$ is {\em Gorenstein$^*$} if it is bounded, graded and finite, and every open interval in $P$ has the homology of a sphere of top dimension.
\end{definition}

It follows from the correspondence between the M\"obius function and the Euler characteristic (the theorem of Ph.\ Hall) that the Gorenstein$^*$ property implies the Eulerian property.

\begin{definition}
A simplicial complex is a {\em piecewise linear}, or {\em PL}, {\em sphere} if it admits a subdivision which is a subdivision of the boundary of a simplex.
\end{definition}
\noindent In particular, PL spheres are of course homeomorphic to spheres. The former class is more well-behaved under certain operations, and this sometimes facilitates inductive arguments. Our interest in PL spheres, rather than spheres in general, will be confined to such situations.

A {\em simplicial poset} is a finite poset equipped with $\Hn$ in which every interval is isomorphic to a Boolean lattice. Such a poset is the {\em face poset} (i.e.\ poset of cells ordered by inclusion) of a certain kind of regular CW complex called {\em Boolean cell complex}. Thus, Boolean cell complexes are slightly more general than simplicial ones, since we allow simplices to share vertex sets. Such complexes were first considered by Bj\"orner \cite{bjorner} and by Garsia and Stanton \cite{GS}.

A Boolean cell complex is {\em pure} if its {\em facets} (inclusion-maximal cells) are equidimensional. A pure complex is {\em thin} if every cell of codimension $1$ is contained in exactly two facets.

A {\em shelling order} of a pure Boolean cell complex $\Delta$ is an ordering $F_0, \dots, F_t$ of the facets of $\Delta$ such that the subcomplex $F_i \cap (\cup_{\alpha < i}F_\alpha)$ is pure of codimension $1$ for all $i \in [t]=\{1,\dots,t\}$. If $\Delta$ has a shelling order, then it is {\em shellable}. The next result can be found in \cite{bjorner}.
\begin{proposition}\label{pr:thin}
If $\Delta$ is a pure, thin, finite and shellable Boolean cell complex, then $\Delta$ is homeomorphic to a sphere. 
\end{proposition}

Let $f_i$ be the number of $i$-dimensional cells in $\Delta$ (including $f_{-1} = 1)$. The {\em $f$-polynomial} is the generating function of the $f_i$:
\[
f_\Delta(x) = \sum_{i\geq 0}f_{i-1}x^i.
\] 
An equivalent, often more convenient, way of encoding this information is the {\em $h$-polynomial}. Setting $d = \dim(\Delta)$, it is defined by
\[
h_\Delta(x) = (1-x)^{d+1}f_\Delta\left(\frac{x}{1-x}\right).
\]
We define coefficients $h_i$ by 
\[
h_\Delta(x) = \sum_{i=0}^{d+1}h_ix^i.
\]
The next result follows from \cite[Proposition 4.4]{stanley}.
\begin{proposition}[Dehn-Sommerville equations]
Suppose $\Delta$ is a pure Boolean cell complex of dimension $d$ which is homeomorphic to a sphere. Then, for all $i\in [d+1]$,
\[
h_{d+1-i} = h_i.
\]
\end{proposition}

If $\Delta$ is shellable, the $h_i$ have nice combinatorial interpretations that we now describe. Suppose $F_1, \dots, F_t$ is a shelling order of $\Delta$. It can be proved that for all $i\in [t]$, $\cup_{\alpha \leq i}F_i$ contains a unique minimal cell which is not contained in $\cup_{\alpha < i}F_\alpha$. Let $r_i$ denote the dimension of this cell. Then, for all $j$,
\[
h_j = |\{i\in[t]\mid r_i = j-1\}|.
\]

\subsection{Properties of Coxeter groups}
We now briefly review some important features of Coxeter groups for later use. We refer the reader to \cite{BB} or \cite{humphreys} for material which may not be familiar.

Henceforth, let $(W,S)$ be a Coxeter system with $|S|<\infty$. The Coxeter length function is $\ell:W\to \N$. If $w = s_1\dots s_k \in W$ and $\ell(w) = k$, the word $s_1\dots s_k$ is called a {\em reduced expression} for $w$. Here and in what follows, symbols of the form $s_i$ are always assumed to be elements in $S$. We do not distinguish notationally between a word in the free monoid over $S$ and the element in $W$ that it represents; we trust the context to make the meaning clear.

Given $w\in W$, we define the {\em left} and {\em right descent sets}, respectively, by
\[
D_L(w) = \{s\in S\mid \ell(sw)<\ell(w)\}
\]
and
\[
D_R(w) = \{s\in S\mid \ell(ws)<\ell(w)\}. 
\]
Observe that $s\in D_R(w)$ (resp.\ $D_L(w)$) iff $w$ has a reduced expression ending (resp.\ beginning) with $s$.

The following two results are equivalent formulations of an important structural property of Coxeter groups. The first has an obvious analogous formulation for left descents.
\begin{proposition}[Exchange property]\label{pr:exchange1}
Suppose $s_1\dots s_k$ is a reduced expression for $w\in W$. If $s\in D_R(w)$, then $ws = s_1\dots\wh{s_i}\dots s_k$ for some $i\in [k]$ (where the hat denotes omission).
\end{proposition}
\begin{proposition}[Deletion property]
If $w = s_1\dots s_k$ and $\ell(w) < k$, then $w = s_1\dots\wh{s_i}\dots\wh{s_j}\dots s_k$ for some $1\leq i<j\leq k$.
\end{proposition}

The most important partial order on $W$ is probably the one we now define. Denote by $T$ the set of {\em reflections} in $W$, i.e.\
\[
T = \{wsw^{-1}\mid w\in W\text{ and }s\in S\}.
\]
\begin{definition}\label{de:bruhat}
The {\em Bruhat order} on $W$ is the partial order $\leq$ defined by $v\leq w$ iff there exist $t_1,\dots, t_k\in T$ such that $w = vt_1\dots t_k$ and $\ell(vt_1\dots t_{i-1}) < \ell(vt_1\dots t_i)$ for all $i\in [k]$. We denote this poset by $\Br(W)$.
\end{definition}

Although not immediately obvious from Definition \ref{de:bruhat}, the Bruhat order is graded with rank function $\ell$. It is topologically well-behaved:
\begin{theorem}[\cite{bjorner, BW, dyer}]
The open intervals in $\Br(W)$ are PL spheres.
\end{theorem}

Innocent as it seems, the next property is nevertheless the key to many results on $\Br(W)$. It follows from Deodhar \cite[Theorem 1.1]{deodhar}. Again, there is an obvious formulation for left descents.
\begin{proposition}[Lifting property]\label{pr:lifting}
Let $v,w\in W$ with $v\leq w$ and suppose $s\in D_R(w)$. Then,
\begin{enumerate}
\item[(i)] $vs \leq w$.
\item[(ii)] $s\in D_R(v) \Rightarrow vs\leq ws$.
\end{enumerate}
\end{proposition}

Next, we define another fruitful way to order $W$. It is readily seen that the following is a weaker order than $\Br(W)$:
\begin{definition}
The {\em right weak order} on $W$ is the partial order $\wkW$ defined by $v\wkW w$ iff $w = vu$ for some $u\in W$ with $\ell(u) = \ell(w) -\ell(v)$. 
\end{definition}
\noindent There is of course also a left weak order defined in the obvious way; we denote it by $\Wwk$. Clearly, both weak orders are graded with rank function $\ell$. 

\section{The combinatorics of twisted involutions}\label{se:comb}
Let $(W,S)$ be a Coxeter system, and $\theta: W \to W$ a group automorphism such that $\theta^2 = \id$ and $\theta(S)=S$. In other words, $\theta$ is induced by an involutive automorphism of the Coxeter graph of $W$.
\begin{definition}
The set of {\em twisted involutions} is
\[
\twist(\theta) = \{w \in W\mid \theta(w) = w^{-1}\}.
\]
\end{definition}
\noindent Note that $\twist(\id)$ is the set of ordinary involutions in $W$.

In this section, we show that the combinatorics of $\twist(\theta)$ is strikingly similar to that of $W$. Our results in this section have been developed for finite $W$ by Springer \cite{springer} and by Richardson and Springer \cite{RS,RS2}, but their proofs (specifically, of the crucial \cite[Lemma 8.1]{RS}) do not hold in the general case, since they make use of the existence of a longest element in $W$. 

\begin{example}[cf.\ Example 10.1 in \cite{RS}]\label{ex:generalization}
{\em Let $W$ be any Coxeter group, and consider the automorphism $\theta:W\times W\to W\times W$ given by $(v,w)\mapsto (w,v)$. It is easily seen that
\[
\twist(\theta) = \{(w,w^{-1})\mid w\in W\},
\]
so that we have a natural bijection $\twist(\theta) \longleftrightarrow W$. This construction allows many properties of $\twist(\theta)$ to be seen as generalizations of Coxeter group properties.}
\end{example}

Consider the set of symbols $\SSS = \{\s \mid s \in S\}$. The free monoid over $\SSS$ acts from the right on the set $W$ by 
\[
w\s = 
\begin{cases}
ws & \text{if $\theta(s)ws = w$,}\\
\theta(s)ws & \text{otherwise},
\end{cases}
\]
and $w\s_1\dots \s_k = (\dots((w\s_1)\s_2)\dots)\s_k$. Observe that $w\s\s=w$ for all $w\in W$, $s\in S$. By abuse of notation, we write $\s_1\dots \s_k$ instead of $e\s_1\dots \s_k$, where $e\in W$ is the identity element. 

\begin{remark}
{\em In \cite{RS}, a slightly different monoid action is used. It satisfies the relation $\s\s = \s$ rather than $\s\s = 1$. We use our formulation since it makes the results easier to state and the similarity to the situation in Coxeter groups more transparent. }
\end{remark}

In $\twist(\theta)$ it is sometimes more convenient to use the following equivalent definition of the action:
\begin{lemma}\label{le:altdef}
Suppose $w\in \twist(\theta)$ and $s\in S$. Then,
\[
w\s = 
\begin{cases}
ws & \text{if $\ell(\theta(s)ws) = \ell(w)$,}\\
\theta(s)ws & \text{otherwise}.
\end{cases}
\]
\begin{proof}
Clearly, $w\s = ws$ implies $\ell(\theta(s)ws) = \ell(w)$. Conversely, suppose that $\ell(\theta(s)ws) = \ell(w)$. If $\theta(s)ws \neq w$, $w$ must have a reduced expression which begins with $\theta(s)$ or ends with $s$. Assume without loss of generality that $\theta(s)s_1\dots s_k$ is such an expression. Since $\theta(w)=w^{-1}$, we have $\ell(ws)< \ell(w)$. No reduced expression for $w$ can both begin with $\theta(s)$ and end with $s$; the Exchange property therefore implies $ws = s_1\dots s_k$, so that $\theta(s)ws = w$. 
\end{proof}
\end{lemma}

Our interest in this action stems from the fact that the orbit of the identity element is precisely $\twist(\theta)$, as the next proposition shows.

\begin{proposition}~

\begin{enumerate}
\item[(i)] For all $\s_1, \dots, \s_k\in \SSS$, we have $\s_1\dots \s_k \in \twist(\theta)$.
\item[(ii)] Given $w \in \twist(\theta)$, there exist symbols $\s_1, \dots, \s_k \in \SSS$ such that $w = \s_1 \dots \s_k$.
\end{enumerate}
\begin{proof}
It is readily checked that $w\s \in \twist(\theta)$ iff $w \in \twist(\theta)$. This proves (i). Noting that $\ell(ws) < \ell(w) \Leftrightarrow \ell(w\s) < \ell(w)$, (ii) follows by induction over the length.
\end{proof}
\end{proposition}

Motivated by this proposition, we define the {\em rank} $\rho(w)$ of a twisted involution $w\in \twist(\theta)$ to be the minimal $k$ such that $w = \s_1\dots \s_k$ for some $\s_1, \dots, \s_k \in \SSS$. The expression $\s_1\dots\s_k$ is then called a {\em reduced $\SSS$-expression} for $w$.

\begin{example}
{\em Let $W$ and $\theta$ be as in Example \ref{ex:generalization}. If $s_1\dots s_k$ is a reduced expression for $w\in W$, then $\underline{(s_1,e)}\dots\underline{(s_k,e)}$ is a reduced $\SSS$-expression for the corresponding twisted involution $(w,w^{-1})\in \twist(\theta)$. To see that it is reduced, note that by construction $\rho$ is always at least half the length, and that the length of $(w,w^{-1})$ in $W\times W$ is $2k$.}
\end{example}

We write $\Br(\twist(\theta))$ for the subposet of $\Br(W)$ induced by $\twist(\theta)$. The study of $\Br(\twist(\theta))$ was initiated in \cite{RS,RS2} because of its connection to Bruhat decompositions of certain symmetric varieties. 

It follows from \cite[Theorem 4.8]{hultman} that $\Br(\twist(\theta))$ is graded with rank function $\rho$ and that $\rho(w) = (\ell(w) + \ell^\theta(w))/2$ for all $w\in \twist(\theta)$, where $\ell^\theta$ is the twisted absolute length function (see \cite{hultman} for the definition). When $W$ is finite, a different way to define the rank function was provided by Richardson and Springer \cite{RS}; the equivalence between the two formulations is due to Carter \cite[Lemma 2]{carter}. In the case of $W$ being a classical Weyl group and $\theta = \id$, Incitti \cite{incitti, incitti2, incitti3} found the rank function using combinatorial arguments.

\begin{example}
{\em With $W$ and $\theta$ as in Example \ref{ex:generalization}, it is readily seen that $\Br(\twist(\theta)) \cong \Br(W)$. Thus, ordinary Bruhat orders are a special case of this construction.}
\end{example}

We now proceed to prove a number of facts that are completely analogous to familiar results from the theory of Coxeter groups.

\begin{lemma}\label{le:rank}
For all $w\in \twist(\theta)$, $s\in S$, we have $\rho(w\s)=\rho(w)\pm 1$, and $\rho(w\s) = \rho(w)-1 \Leftrightarrow s\in D_R(w)$.
\begin{proof}
Since $w\s\s = w$, $|\rho(w\s)-\rho(w)|\leq 1$. Using Lemma \ref{le:altdef}, it is straightforward to verify from Definition \ref{de:bruhat} that $w > w\s$ iff $s\in D_R(w)$, and otherwise $w<w\s$. The lemma now follows from the fact that $\rho$ is the rank function of $\Br(\twist(\theta))$.
\end{proof}
\end{lemma}

\begin{lemma}[Lifting property for $\SSS$] \label{le:lift}
Let $s\in S$ and $v,w\in W$ with $v \leq w$, and suppose $s\in D_R(w)$. Then,
\begin{enumerate}
\item[(i)] $v\s \leq w$.
\item[(ii)] $s\in D_R(v) \Rightarrow v\s \leq w\s$.
\end{enumerate}
\begin{proof}
Suppose $s\in D_R(v)$. If $v\s = vs$ and $w\s = \theta(s)ws$, Proposition \ref{pr:lifting} shows first that $vs \leq ws$, then $v = \theta(s)vs \leq ws$ and, finally, $v\s = \theta(s)v \leq \theta(s)ws = w\s$. The other cases admit similar proofs. 
\end{proof}
\end{lemma}

\begin{proposition}[Exchange property for $\twist(\theta)$]\label{pr:exchange2}
Suppose $\s_1 \dots \s_k$ is a reduced $\SSS$-expression and that $\rho(\s_1 \dots \s_k \s) < k$. Then, $\s_1 \dots \s_k \s = \s_1\dots \wh{\s_i}\dots \s_k$ for some $i \in [k]$.
\begin{proof}
Let $w = \s_1 \dots \s_k$ and $v = \s_1 \dots \s_k \s$. We have $w > v$. Let $i \in [k]$ be maximal such that $v\s_k\dots\s_i>v\s_k\dots\s_{i+1}$ (it exists; otherwise we would have $\rho(v\s_k\dots\s_1)<0$). Repeated application of Lemma \ref{le:lift} shows that $w\s_k\dots\s_{i+1}\geq v\s_k\dots\s_i$. Since $\rho(w\s_k\dots\s_{i+1}) = \rho(v\s_k\dots\s_i)$ and $\rho$ is the rank function of $\Br(\twist(\theta))$, this implies $w\s_k\dots\s_{i+1} = v\s_k\dots\s_i$. Thus, $v =  w\s_k\dots\s_{i+1}\s_i\dots\s_k =  \s_1\dots \wh{\s_i}\dots \s_k$.
\end{proof}
\end{proposition}

\begin{proposition}[Deletion property for $\twist(\theta)$]\label{pr:deletion2}
If $\rho(\s_1\dots \s_k) < k$, then we have $\s_1 \dots \s_k = \s_1\dots \wh{\s_i}\dots \wh{\s_j} \dots \s_k$ for some $1\leq i < j \leq k$.
\begin{proof}
Let $j\in [k]$ be minimal such that $\s_1\dots\s_j$ is not reduced and apply Proposition \ref{pr:exchange2} to this expression.
\end{proof}
\end{proposition}

\section{The Bruhat order on twisted involutions}\label{se:bruhat}
We now turn our attention to the poset $\Br(\twist(\theta))$. Incitti \cite{incitti, incitti2, incitti3} showed that $\Br(\twist(\id))$ is EL-shellable and Eulerian whenever $W$ is a classical Weyl group. He conjectured the same to hold for any Coxeter group $W$ (when $W$ is infinite, it should hold for every interval in $\Br(\twist(\id))$). In \cite{hultman}, it was proved that every interval in $\Br(\twist(\theta))$ is Gorenstein$^*$ over $\Z_2$ (in particular Eulerian) for arbitrary $W$ and $\theta$. The purpose of this section is to strengthen this result by showing that every interval in $\Br(\twist(\theta))$ is a PL sphere. Thus, the main consequence of Incitti's conjecture holds, although it is still not proved that the spheres actually are shellable. 

A key poset property is that of admitting a {\em special matching}:
\begin{definition}
Let $P$ be a poset and $M:P\to P$ an involution such that for all $x\in P$, either $x$ covers $M(x)$ or $M(x)$ covers $x$. Then, $M$ is called a {\em special matching} iff for all $x,y\in P$ such that $M(x)\neq y$, it holds that
\[
y\text{ covers }x \Rightarrow M(x) < M(y).
\]
\end{definition}
The term ``special matching'' is due to Brenti, see \cite{brenti}. In Eulerian posets, special matchings are equivalent to {\em compression labellings} in the sense of du Cloux \cite{ducloux}. Although the statement of \cite[Corollary 3.6]{ducloux} is slightly weaker, its proof implies the next result. It is an application of \cite[Theorem 3.5]{ducloux}, which, in turn, is a reformulation of a result from Dyer's thesis \cite{dyer}. It was reproved in the setting of Bruhat orders by Reading \cite{reading}, and his proof is easily adapted to the general situation.
\begin{theorem}[\cite{ducloux,dyer,reading}]\label{th:sphere}
Let $P$ be an Eulerian poset with a special matching $M$. If $(\Hn,M(\He))$ is a PL sphere, then so is $\oP$.
\end{theorem}
\begin{corollary}\label{co:zircon}
Suppose $P$ is an Eulerian poset in which every lower interval $[\Hn,x]$, $x \neq \Hn$, has a special matching. Then, every open interval in $P$ is a PL sphere.
\begin{proof}
Links in PL spheres are PL spheres (see \cite[Theorem 4.7.21.iv]{BLSWZ}). Therefore, it is enough to prove that $\oP$ is a PL sphere. This follows from Theorem \ref{th:sphere} by induction over the rank.
\end{proof}
\end{corollary}

\begin{remark}{\em
Explicitly requiring $P$ to be Eulerian in Corollary \ref{co:zircon} is not important. In fact, if every lower interval in a bounded poset $P$ has a special matching, then $P$ is necessarily Eulerian \cite{BCM}.}
\end{remark}

\begin{theorem}\label{th:main}
Let $w\in \twist(\theta)$, and suppose $s\in D_R(w)$. Then, the map $v \mapsto v\s$ is a special matching on the interval $[e,w]\subseteq \Br(\twist(\theta))$.
\begin{proof}
Part (i) of Lemma \ref{le:lift} shows that $v\mapsto v\s$ maps $[e,w]$ to itself, and the fact that $v\s\s = v$ for all $v$ shows that the map is an involution. As in the proof of Lemma \ref{le:rank}, $v$ and $v\s$ are always comparable. By rank considerations, one of them must therefore cover the other.

Now, pick $x,y\in [e,w]$ such that $y$ covers $x$ and consider Lemma \ref{le:lift}. If $x\s < x$, part (ii) shows that $x\s < y\s$. If $x\s> x$ and $y\s < y$, we must have $x\s = y$ by part (i). Finally, if $x\s > x$ and $y\s > y$, then $x\s < y\s$, again by part (i). Thus, $v \mapsto v\s$ is a special matching. 
\end{proof}
\end{theorem}

\begin{corollary}\label{co:PL}
The open intervals in $\Br(\twist(\theta))$ are PL spheres.
\end{corollary}

\begin{remark}{\em Using \cite[Proposition 4.7.23]{BLSWZ}, we may deduce from Corollary \ref{co:PL} that every half-open interval $[v,w)\subseteq \Br(\twist(\theta))$ is the face poset of a regular cell decomposition of a sphere (in which the empty set, corresponding to $v$, is regarded to be a cell).}
\end{remark}

It follows from Theorem \ref{th:main} that the intervals in $\Br(\twist(\theta))$ are {\em accessible} posets as defined in \cite{ducloux}. It is known that not all accessible posets are Bruhat intervals, i.e.\ intervals in some $\Br(W)$. Interestingly, the smallest counterexamples (called {\em Dyer obstructions} in \cite{ducloux}) coincide with $\Br(\twist(\theta))$ when $W=A_3$ and $\theta$ is either the identity or the unique non-trivial Coxeter graph automorphism, respectively. This makes it natural to ask whether or not all accessible posets arise as intervals in Bruhat orders on twisted involutions. 

\section{A Coxeter complex analogue for twisted involutions}\label{se:weak}
In this section, we construct a cell complex $\Delta_\theta$ whose relationship with $\twist(\theta)$ has many features in common with the connection between the Coxeter complex $\Delta_W$ and $W$. Although some results also make sense for infinite groups, our main interest here is in the finite setting. Therefore, throughout the rest of the paper, $(W,S)$ will be a {\em finite} Coxeter system with an involutive automorphism $\theta$. 

We define a graph $G_\theta$ on the vertex set $\twist(\theta)$ with edges labelled by elements in $S$ as follows: there is an edge with label $s$ between $v$ and $w$ iff $v\s = w$. 

If we direct all edges according to decreasing $\rho$-values and merge multiple edges, we obtain from $G_\theta$ the Hasse diagram of the following partial order which was first defined in \cite{RS}: 

\begin{definition}
The {\em weak order} on $\twist(\theta)$ is the partial order $\wk$ defined by $v \wk w$ iff there exist  $\s_1, \dots, \s_k\in \SSS$ such that $v\s_1\dots\s_k = w$ and $\rho(v) = \rho(w) - k$. We denote this poset by $\Wk(\theta)$. 
\end{definition}

Observe that $\Wk(\theta)$ is a subposet of $\Br(\twist(\theta))$, i.e.\ the identity map $\Wk(\theta) \to \Br(\twist(\theta))$ is order-preserving.

It should be noted that $\Wk(\theta)$ does not in general coincide with the order induced by the (left or right) weak order on $W$. In fact, while the former is clearly graded with rank function $\rho$, the latter is not a graded poset in general.

\begin{example}\label{ex:leftright}
{\em Return to the situation in Example \ref{ex:generalization}. Observe that 
\[
(w,w^{-1})\underline{(s,e)} = (ws,(ws)^{-1}) \wk (w,w^{-1}) \Leftrightarrow s\in D_R(w) \Leftrightarrow ws \wkW w
\]
and
\[
(w,w^{-1})\underline{(e,s)} = (sw,(sw)^{-1}) \wk (w,w^{-1}) \Leftrightarrow s\in D_L(w) \Leftrightarrow sw \Wwk w.
\]
Hence, in this setting, $\Wk(\theta)$ is isomorphic to the transitive closure of the union of the left and right weak orders on $W$. This poset is sometimes called the {\em two-sided weak order} on $W$. It was studied by Bj\"orner in \cite{bjorner3}.
}
\end{example}

Given $J\subseteq S$, consider the subgraph of $G_\theta$ obtained by removing all edges with labels not in $J$. For $w \in \twist(\theta)$, let $wC_J$ be the connected component which contains $w$ in this subgraph. It should be stressed that we regard $wC_J$ as an edge-labelled graph, not merely as a set of vertices. Define
\[
P_\theta = \{wC_J\mid w\in \twist(\theta)\text{ and }J\subseteq S\}.
\]
The elements of $P_\theta$ are partially ordered by reverse inclusion, i.e.\ $g_1 \leq g_2$ iff $g_2$ is a (labelled) subgraph of $g_1$.
\begin{proposition}
The poset $P_\theta$ is the face poset of a pure Boolean cell complex $\Delta_\theta$ of dimension $|S|-1$.
\begin{proof}
The bottom element of $P_\theta$ is $G_\theta$. The maximal elements are the twisted involutions. Let $w\in \twist(\theta)$. The map $wC_J \mapsto J$ is easily seen to be a poset isomorphism from the interval $[G_\theta, w] = [wC_S,wC_\emptyset]\subseteq P_\theta$ to the dual of the Boolean lattice of subsets of $S$.
\end{proof}
\end{proposition}

\begin{remark}
{\em We briefly indicate why we regard $\Delta_\theta$ as a Coxeter complex analogue. Suppose we replace $\twist(\theta)$ with $W$ and $G_\theta$ with the Cayley graph of $W$ (with respect to the generating set $S$). The connected component $wC_J$ would then become the subgraph induced by the parabolic coset $w\langle J\rangle$. Ordering the set of such cosets by reverse inclusion produces the face poset of the Coxeter complex $\Delta_W$. We refer to Brown's book \cite{brown} for a thorough background on Coxeter complexes.}
\end{remark}

\begin{remark}
{\em The complex $\Delta_\theta$ is not in general a simplicial complex. For example, if $s$ and $\theta(s)\neq s$ commute for some $s\in S$, we have $\s = \underline{\theta(s)} = s\theta(s)$. Thus, there are two edges between $e$ and $s\theta(s)$ in $G_\theta$, implying that the facets in $\Delta_\theta$ indexed by $e$ and $s\theta(s)$ share two codimension $1$ cells. Similar examples exist when $\theta =\id$.}
\end{remark}

\begin{lemma}\label{le:minimal}
Let $J\subseteq S$ and $w\in W$. Then, the vertex set of $wC_J$ contains a unique $\Wk(\theta)$-minimal element $\min(w,J)$. 
\begin{proof}
Given $s\in J$, a $\Wk(\theta)$-minimal element in $wC_J$ must clearly not have a reduced $\SSS$-expression ending in $\s$. Conversely, if $v\in wC_J$ is not $\Wk(\theta)$-minimal, there exists $s\in J$ such that $\rho(v\s) < \rho(v)$. Since $v\s\s = v$, we obtain a reduced $\SSS$-expression for $v$ ending in $\s$ by attaching $\s$ to any reduced $\SSS$-expression for $v\s$. Thus, the $\Wk(\theta)$-minimal elements in $wC_J$ are precisely the elements that have no reduced $\SSS$-expressions ending in $\s$ for all $s\in J$. It is clear that at least one such element exists.

Now suppose $u$ and $v$ are two $\Wk(\theta)$-minimal elements in $wC_J$. Let $\s_1\dots\s_k$ be a reduced $\SSS$-expression for $v$. We may write
\[
u = v\s_{k+1}\dots\s_{k+l} =\s_1\dots\s_{k+l},
\]
for some $s_{k+1},\dots,s_{k+l}\in J$. By Proposition \ref{pr:deletion2}, this expression contains a reduced subexpression for $u$. Since it cannot end in any $\s_t$, $k+1\leq t\leq k+l$, it must be a subexpression of $\s_1\dots\s_k$. By symmetry, this subexpression must, in turn, contain a reduced subexpression for $v$. Thus, $u = v$.
\end{proof}
\end{lemma}

The completely analogous Coxeter complex version of the next result is due to Bj\"orner \cite[Theorem 2.1]{bjorner2}.

\begin{theorem}\label{th:shelling}
Any linear extension of $\Wk(\theta)$ is a shelling order for $\Delta_\theta$. In particular, $\Delta_\theta$ is shellable.
\begin{proof}
Suppose $w_1, \dots, w_k$ (where $k=|\twist(\theta)|$) is a linear extension of $\Wk(\theta)$. Let $j \in [k]$, and suppose $g$ is a cell in $w_i \cap w_j$ for some $i < j$. We must show that $g$ is contained in a codimension 1 cell in $w_j \cap (\cup_{\alpha < j}w_\alpha$). In terms of the graphs that represent the cells, the situation is this: $g=w_iC_J = w_jC_J$ for some $J\subseteq S$. It must be shown that $g$ contains an edge connecting $w_j$ with some $w_\alpha$, $\alpha < j$. If $J$ contains a right descent $s$ of $w_j$, we can use $w_\alpha = w_j\s$. Otherwise, $w_j = \min(w,J)$. By Lemma \ref{le:minimal}, this implies $w_j \wk w_i$, contradicting our choice of linear extension.   
\end{proof}
\end{theorem}
\begin{corollary}
The complex $\Delta_\theta$ is homeomorphic to the $(|S|-1)$-dimensional sphere.
\begin{proof}
Since codimension $1$ cells in $\Delta_\theta$ correspond to edges in $G_\theta$, $\Delta_\theta$ is thin. The corollary now follows from Proposition \ref{pr:thin}. 
\end{proof}
\end{corollary}
We now define an analogue of the $W$-Eulerian polynomial (i.e.\ the generating function counting the elements of $W$ with respect to the number of descents) for twisted involutions:
\[
\des_\theta(x) = \sum_{w \in \twist(\theta)} x^{|D_R(w)|}.
\]
\begin{example}
{\em Again, consider the setting of Example \ref{ex:generalization}. Let $\phi$ denote the natural bijection $W \to \twist(\theta)$. From the argument in Example \ref{ex:leftright}, it follows that $|D_R(\phi(w))| = |D_R(w)| + |D_L(w)|$ for all $w\in W$. Thus, we have in this setting
\[
\des_\theta(x) = \sum_{w\in W} x^{|D_R(w)|+|D_L(w)|} = \sum_{w\in W} x^{|D_R(w)|+|D_R(w^{-1})|}.
\] 
This can be viewed as a two-sided analogue of the $W$-Eulerian polynomial.}
\end{example}

For Coxeter complexes, the counterpart of the next result is \cite[Theorem 2.3]{brenti2}.
\begin{theorem}\label{th:h-vector}
The $h$-polynomial of $\Delta_\theta$ coincides with $des_\theta(x)$ .
\begin{proof}
Consider the shelling order of $\Delta_\theta$ given in the proof of Theorem \ref{th:shelling}. The unique minimal new cell introduced in the $i$th shelling step is $w_iC_J$, where $J$ is the set of right ascents (i.e.\ non-descents) of $w_i$. The dimension of this cell is $|S|-|J|-1 = |D_R(w_i)|-1$. For the $h$-vector of $\Delta_\theta$, this means that
\[
h_j = |\{w\in \twist(\theta)\mid |D_R(w)|=j\}|.
\]
Thus, 
\[
\des_\theta(x) = \sum_{j=0}^{|S|} h_j x^j.
\]
\end{proof}
\end{theorem}
A polynomial $P\in \Z[x]$ is called {\em symmetric} if $x^dP(x^{-1}) = P(x)$, where $d = \deg(P)$.
\begin{corollary}\label{co:symmetric}
The polynomial $\des_\theta(x)$ is symmetric.
\begin{proof}
This is immediate from Theorem \ref{th:h-vector} and the Dehn-Sommerville equations.
\end{proof}
\end{corollary}
\begin{remark}
{\em Suppose $W$ is irreducible, and let $w_0$ denote the longest element in $W$. It is known (see \cite[Exercise 4.10]{BB}) that $ww_0 = w_0w$ for all $w\in W$ unless $W$ is of one of the types $I_2(2n+1), A_n$, $D_{2n+1}$ and $E_6$. Thus, in all other cases $w\mapsto ww_0$ is an involution $\twist(\theta)\to \twist(\theta)$ which sends ascents to descents, proving Corollary \ref{co:symmetric} for these cases. When $W = A_n$, $\theta=\id$, Corollary \ref{co:symmetric} is due to Strehl \cite{strehl}. See also \cite{dukes}.}
\end{remark}

\end{document}